 
\input amstex
\magnification 1200
\loadmsbm
\parindent 0 cm

\define\nl{\bigskip\item{}}
\define\snl{\smallskip\item{}}
\define\inspr #1{\parindent=20pt\bigskip\bf\item{#1}}
\define\iinspr #1{\parindent=27pt\bigskip\bf\item{#1}}
\define\einspr{\parindent=0cm\bigskip}


\input amssym
\input amssym.def

\centerline{\bf The Fourier transform in quantum group theory \rm($^*$)}
\bigskip\bigskip
\centerline{\it A.\ Van Daele \rm ($^{**}$)}
\bigskip\bigskip\bigskip
{\bf Abstract} 
\bigskip 
The Fourier transform, known in classical analysis, and generalized in abstract harmonic analysis, can also be considered in the theory of locally compact quantum groups.
\snl
In this talk, I will discuss some aspects of this more general Fourier transform. In order to avoid technical difficulties, typical for the analytical approach, I will restrict to the algebraic quantum groups. Roughly speaking, these are  the locally compact quantum groups that can be treated with purely algebraic methods (in the framework of multiplier Hopf algebras). 
\snl
I will illustrate various notions and results using not only classical Fourier theory on the circle group $\Bbb T$, but also on the additive group $\Bbb Q_p$ of $p$-adic numbers. It should be observed however that these cases are still too simple to illustrate the full power of the more general theory.
\nl
\nl
{\it April 2007} (Version 1.1)
\vskip 6 cm
\hrule
\bigskip\parindent 0.8 cm
\item{($^*$)} Talk at the meeting {\it New techniques in Hopf algebras and graded ring theory}. Brussels, September 19-23, 2006.
\item{($^{**}$)} Department of Mathematics, K.U.\ Leuven, Celestijnenlaan 200B,
B-3001 Heverlee (Belgium). E-mail: Alfons.VanDaele\@wis.kuleuven.be
\parindent 0 cm

\newpage

\bf 0. Introduction \rm
\nl
Let $G$ be an abelian locally compact group. The dual group $\widehat G$ is defined as the group of continuous characters on $G$ (with values in the circle group $\Bbb T=\{ z\in \Bbb C \mid |z|=1\}$). With a suitable topology, $\widehat G$ is again a locally compact abelian group. Pontryagin's duality states that the dual of $\widehat G$ is again $G$. Typical examples are the groups $\Bbb Z$, $\Bbb T$ and $\Bbb R$ with duals $\Bbb T$, $\Bbb Z$ and $\Bbb R$ respectively. We are also interested in the additive group $\Bbb Q_p$ of $p$-adic numbers whose dual is again $\Bbb Q_p$. 
\snl
The Fourier transform of a continuous complex function $f$ with compact support on such a group $G$ is a continuous complex function $\widehat f$ on $\widehat G$, defined by 
$$\widehat f(q)=\int_G f(r)\langle r,q\rangle^-\,dr$$
where integration is over the Haar measure on $G$, where $\langle r,q\rangle$ is used to denote the value of the character $q$ in the point $r$ and where $z^-$ stands for the complex conjugate of the complex number $z$. In general, the function $\widehat f$ will no longer have compact support, but it will tend to $0$ at infinity. In the case $G=\Bbb T$, the Fourier transform is a function on $\Bbb Z$. In the case $G=\Bbb R$, the Fourier transform is again a function on $\Bbb R$. In the first case, we obtain the Fourier coefficients. In the second case, we get the classical Fourier integral. In general, we end up doing abstract harmonic analysis.
\snl
Among the various interesting properties of the Fourier transform, we have the Fourier inversion formula, Plancherel's formula and the fact that convolution product is converted to pointwise product by the Fourier transform.
\snl
See e.g.\ [F]  for a nice and comprehensive treatment of the theory of abstract harmonic analysis.
\nl
If the group $G$ is no longer assumed to be abelian, it is not possible anymore to consider the dual group $\widehat G$. For a long time, people have tried to construct objects in order to generalize Pontryagin's theorem to the non-abelian case. Recently, this problem found a nice and satisfactory solution with the introduction of the theory of locally compact quantum groups. It is possible to construct a dual of the same type and in this way, the theorem of Pontryagin can be generalized in a beautiful way to the non-abelian case. Also the Fourier transform has a natural generalization to locally compact quantum groups and many of the classical results can be extended to this more general framework.
\snl
We refer to [K-V1] and [K-V2] for the theory of locally compact quantum groups, see also [VD3].
\nl
This analytical theory is highly non-trivial and technically quite difficult. It uses, among other things, a theory of non-commutative integration. This would be far beyond the scope of this talk (and this conference). Fortunately however, there is the purely algebraic theory of multiplier Hopf algebras with integrals (the so-called algebraic quantum groups), having many of the features of the general theory. They can be treated without going into the analytical technical difficulties.
\snl
So, in this talk, I will discuss the concept of the Fourier transform for multiplier Hopf algebras with integrals.  To do Fourier analysis, one needs an integral on $A$ as well as on the dual $\widehat A$. This is the case when $A$ is a multiplier Hopf algebra with integrals. Then, the dual $\widehat A$ can be considered and it is again a multiplier Hopf algebra with integrals. The theory of Hopf algebras is not sufficiently general for this purpose because in this case, requiring an integral both on $A$ and its dual, forces $A$ to be finite-dimensional. Therefore, many interesting cases and examples can not be treated if we stick to the theory of Hopf algebras.
\nl
In {\it Section} 1, we introduce the main concept. First we recall briefly the notion of an algebraic quantum group $(A,\Delta)$ and its dual $(\widehat A,\widehat \Delta)$. We introduce the notion of the Fourier transform in this setting and we motivate our definition within the classical framework of the theory of Fourier series. In {\it Section} 2, we prove some of the main properties of this Fourier transform. We treat the Fourier inversion formula, Plancherel's formula and we show that convolution product on $A$ is converted to the usual product on $\widehat A$. We discuss the notion of algebraic quantum groups of discrete and of compact type and show that these types are dual to each other. Finally, in {\it Section} 3, we focuss on an example of a totally disconnected locally compact group, the additive group $\Bbb Q_p$ of the $p$-adic numbers. We introduce the associated multiplier Hopf algebra. Because it is highly related, we introduce, again in general, the notion of a group-like projection and we show that the Fourier transform of such a group-like projection is again a group-like projection. We illustrate the result in the case of the $p$-adic numbers. We finish with a short section ({\it Section} 4) where we draw some conclusions and discuss related research.
\nl
We have already given some {\it basic references} for the classical theory of the Fourier transform (in abstract harmonic analysis), namely [F]. We have also refered to the original works of Kustermans and Vaes [K-V1] and [K-V2], and to [VD3] for the general theory of locally compact quantum groups. For the theory of multiplier Hopf algebras and algebraic quantum groups, we have the original works [VD1] and [VD2], as well as the survey paper [VD-Z]. Finally, there are the papers by Landstad and myself about group-like projections and related things ([L-VD1] and [L-VD2]).  
\nl\nl
\bf Acknowledgements \rm
I would like to thank the organizers of this meeting for giving me the opportunity to talk about this work. I am also grateful to my coworkers in Leuven for discussions about the example of the $p$-adic numbers.
\nl\nl

\bf 1. The Fourier transform for algebraic quantum groups \rm
\nl
Let us first recall the notion of an algebraic quantum group (a regular multiplier Hopf algebra with integrals).

\inspr{1.1} Definition \rm
Let $A$ be an algebra over $\Bbb C$, with or without identity, but with a non-degenerate product. Let $\Delta$ be a comultiplication on $A$ and assume that it is regular and full. Then the pair $(A,\Delta)$ is called an {\it algebraic quantum group} if there exists faithful left and right integrals on $A$.
\einspr

A coproduct is assumed to satisfy $\Delta(a)(1\otimes b)\in A \otimes A$ and $(a \otimes 1)\Delta(b)\in A\otimes A$ for all $a,b\in A$. If also $\Delta(a)(b\otimes 1)\in A \otimes A$ and $(1 \otimes a)\Delta(b)\in A\otimes A$ for all $a,b\in A$, then $\Delta$ is called {\it regular}. It is called {\it full} if the left and the right legs of $\Delta(A)$ are all of $A$.

\snl
This definition is not the original one, but it refers to the Larson-Sweedler theorem as it is proven for multiplier Hopf algebras in [VD-W]. For the original definition, see [VD1] and [VD2]. Recall that the existence of a counit $\varepsilon$ and an invertible antipode $S$ is proven from these assumptions. They satisfy (and are characterized by) similar equations as in the case of Hopf algebras.

\snl
For any such algebraic quantum group $(A,\Delta)$, we have a dual algebraic quantum group $(\widehat A, \widehat \Delta)$. The result is as follows (see again [VD2]).

\inspr{1.2} Theorem \rm
Let $(A,\Delta)$ be an algebraic quantum group and let $\varphi$ be a left integral. Define $\widehat A$ as the space of linear functionals on $A$ of the form $\varphi(\cdot\,a)$ where $a\in A$. The dual of the coproduct $\Delta$ defines a product on $\widehat A$ and the dual of the product on $A$ defines a coproduct $\widehat\Delta$ on $\widehat A$ making the pair $(\widehat A,\widehat \Delta)$ into an algebraic quantum group. It is called the {\it dual} of $(A,\Delta)$. The dual of $(\widehat A,\widehat \Delta)$ is again $(A,\Delta)$.
\einspr

Recall that integrals are unique (up to a scalar) and so $\widehat A$ is uniquely defined. By the faithfulness of $\varphi$, the space $\widehat A$ separates points of $A$. A right integral $\widehat \psi$ on $\widehat A$ is given by the formula $\widehat\psi(\omega)=\varepsilon(a)$ when $\omega=\varphi(\cdot\,a)$. Similarly, a left integral $\widehat\varphi$ on $\widehat A$ is given by $\widehat\varphi(\omega)=\varepsilon(a)$ if $\omega=\psi(a\,\cdot)$ where $\psi$ is a right integral on $A$. 
\snl
An important case is that of a $^*$-algebra $A$ where $\Delta$ is assumed to be a $^*$-homomorphism. In this case, we call $A$ a $^*$-algebraic quantum group. Then $\widehat A$ is again a $^*$-algebra and the involution is defined by $\omega^*(a)=\omega(S(a)^*)^-$ (where as before, $z^-$ denotes the complex conjugate of the complex number $z$). The dual coproduct $\widehat \Delta$ is again a $^*$-homomorphism. So, the dual is again a $^*$-algebraic quantum group. Later, we will consider such a case where the left integral is assumed to be positive, i.e.\ $\varphi(a^*a)\geq 0$ for all $a$. Often, when we consider $^*$-algebraic quantum groups, we assume positivity of the left integral. 

\nl 
In the setting of algebraic quantum groups, the Fourier transform is defined as follows.

\inspr{1.3} Definition \rm
The {\it Fourier transform} is a linear map $\Cal F$ from $A$ to $\widehat A$, defined by $\Cal F(a)=\varphi(\cdot \, a)$.
\einspr

Sometimes we write $\widehat a$ for $\Cal F(a)$. Remark also that there are many other possible choices one can make to define the Fourier transform. The choice made in Definition 1.3 above is somewhat arbitrary but convenient.

\nl
We will now spend the rest of this section to explain, using the theory with Fourier series, why this terminology is justified.

\inspr{1.4} Example \rm
i) Consider the circle group $\Bbb T$ of complex numbers with modulus $1$. In this case, the Fourier transform maps a continuous complex function $f$ on $\Bbb T$ to a function $\widehat f$ on $\Bbb Z$, defined by
$$\widehat f(n)=\frac{1}{2\pi}\int_0^{2\pi}f(e^{i\theta})e^{-in\theta}\,d\theta$$
where we integrate over the usual Lebesque measure. For all of this to fit into the framework of algebraic quantum groups, we will restrict to the subalgebra $A$ of functions $f$ whose Fourier transform $\widehat f$ has finite support. Such functions are linear combinations of the functions $e_n$, with $n\in \Bbb Z$, defined on $\Bbb T$ by $e_n(z)=z^n$. We do not really loose any relevant information because such functions are uniformly dense in the continuous functions.
\snl
Let us now translate to the framework of algebraic quantum groups.
\snl
ii) The algebra $A$ has an identity and it is made into a Hopf algebra when the coproduct $\Delta$ is defined by 
$$\Delta(f)(z_1,z_2)=f(z_1z_2)$$
for $z_1,z_2\in \Bbb T$. Indeed, $\Delta(e_n)=e_n\otimes e_n$ because $(z_1z_2)^n=z_1^nz_2^n$ for all $z_1,z_2\in \Bbb T$. Of course, $(A,\Delta)$ is nothing else but the Hopf algebra $\Bbb C\Bbb Z$, the group algebra of the additive group $\Bbb Z$ over the complex numbers. A left integral is obtained by integration as above.
\snl
What is the dual of $(A,\Delta)$ in this case?
\snl
iii) Consider the algebra $B$ of complex functions on $\Bbb Z$ with finite support (and pointwise product). This is made into an algebraic quantum group by defining the coproduct in a similar way, namely by
$$\Delta(f)(n,m)=f(n+m)$$
for $f\in B$ and $n,m\in \Bbb Z$. Here, a left integral is obtained by summation over $\Bbb Z$. This is possible because we work with functions with finite support so that the sum over its values exists.
\snl
iv) The pair $(B,\Delta)$ is identified with the dual of $(A,\Delta)$ using the pairing $\langle\,\cdot\, , \,\cdot\,\rangle: A\times B \to \Bbb C$, defined by $\langle e_n,f\rangle =f(-n)$ when $n\in \Bbb Z$ and $f\in B$. Indeed, for $n,m\in \Bbb Z$ and $f,g\in B$, we have
$$\align 
\langle e_n, fg\rangle &= (fg)(-n)=f(-n)g(-n) \\
                       &= \langle e_n \otimes e_n, f\otimes g\rangle \\
                       &= \langle\Delta(e_n),f\otimes g\rangle\\
\langle e_ne_m,f\rangle &= \langle e_{n+m},f\rangle = f(-n-m)\\
                       &= \langle e_n \otimes e_m, \Delta(f)\rangle. 
\endalign$$ 
v) The identification of $\widehat A$  with $B$ is precisely done by means of the Fourier transform. Denote by $\delta_n$ the function on $\Bbb Z$ with value $1$ in $n$ and $0$  everywhere else. Then we have $\Cal F(e_n)= \delta_{n}$ because
$$ \langle e_m,\delta_{n}\rangle = \delta_{n}(-m) = \delta_0(m+n) $$
and also $\varphi(e_me_n)=\varphi(e_{m+n})=\delta_0(m+n)$ for all $n,m\in \Bbb Z$. Remark that, for this identification to be correct, it is important that we have taken the proper joint normalization of the integrals on $A$ and $B$.
\einspr

As one can see, there is really nothing deep going on. It is more a matter of definitions and of notations. The example is too simple (but good enough to illustrate the definition of the Fourier transform above). One should be aware of the fact that, in the non-abelian case, the dual group $\widehat G$ of $G$ does not exist and that the (multiplier) Hopf algebra of functions on $\widehat G$ is replaced by the dual of the (multiplier) Hopf algebra of functions on $G$. This is what illustrates the example above.
\nl\nl

\bf 2. Properties of the Fourier transform \rm
\nl
In this section, we will prove the basic properties of the Fourier transform, as introduced previously.
\snl
Let us first consider the {\it inverse transform}. We have the following result.

\inspr{2.1} Lemma \rm
If $a\in A$ and $\omega=\varphi(\cdot\, a)$, then $\widehat\psi(\omega'\omega)=\omega'(S^{-1}(a))$ for all $\omega'\in \widehat A$.
\snl\bf Proof: \rm
Take $a$ and $\omega$ as above. For any linear functional $\omega'$ in $\widehat A$ and $x\in A$, we have
$$\align (\omega'\omega)(x) &= \omega'(x_{(1)})\omega(x_{(2)}) = \omega'(x_{(1)})\varphi(x_{(2)}a) \\
                           &= \omega'(x_{(1)}a_{(2)}S^{-1}(a_{(1)}))\varphi(x_{(2)}a_{(3)})\\
                           &= \omega'(S^{-1}(a_{(1)}))\varphi(x a_{(2)}).
\endalign$$
Therefore, 
$$\widehat\psi(\omega'\omega)=\varepsilon(a_{(2)})\omega'(S^{-1}(a_{(1)}))=\omega'(S^{-1}(a)).$$
\hfill $\square$
\einspr

We have used the Sweedler notation. Recall that this is justified provided we have the proper 'coverings' (see e.g.\ [VD-Z]).
\snl
This lemma gives the inverse Fourier transform. Indeed, if $\omega=\varphi(\cdot\, a)$, then $a=\widehat\psi(S(\,\cdot\,)\omega)$. We identify, as usual, the dual of $\widehat A$ with $A$ by considering elements in $A$ as linear functionals on $\widehat A$.

\snl
Observe the presence of $S$ here, just as in the classical Fourier inversion formula.

\nl
Let us now consider the {\it convolution product} on $A$ and show that the Fourier transform converts it to multiplication in $\widehat A$.

\inspr{2.2} Proposition \rm
For $a,b\in A$, define the convolution product $a*b\in A$ by
$$a*b=\varphi(S^{-1}(b_{(1)})a)b_{(2)}.$$
Then $\widehat{a*b}=\widehat a \,\widehat b$.
\snl \bf Proof: \rm
For any $x\in A$ we have
$$\align \langle x,\widehat a\, \widehat b \rangle &= (\varphi\otimes\varphi)(\Delta(x)(a\otimes b))\\
                 &= (\varphi\otimes\varphi)(\Delta(xb_{(2)})(S^{-1}(b_{(1)})a\otimes 1)) \\
                 &= \varphi(xb_{(2)}) \varphi(S^{-1}(b_{(1)})a) \\
                 &= \varphi(x(a*b)).
\endalign$$
\hfill $\square$
\einspr

Compare this with the formula for the convolution of continuous functions with compact support on $\Bbb R$. If $f,g$ are two such functions, then their convolution product $f*g$ is again a continuous function with compact support, given by
$$(f*g)(t)=\int_{-\infty}^{+\infty}f(s)g(s^{-1}t)\,ds$$
for all $t\in \Bbb R$.
\snl
Because in general, the algebra $A$ is not abelian and also left and right integrals may be different, there are many possible choices to define this 'convolution product' and for any of these choices, there are several different expressions. For the above choice e.g.\  we can also write
$$a*b=\varphi(S^{-1}(b)a_{(2)})a_{(1)}.$$

\nl
Next, let us show that we have a natural equivalent of {\it Plancherel's formula}. For this, we need to work in the setting of a 
$^*$-algebraic quantum group. Let us assume positivity of $\varphi$.

\inspr{2.3} Proposition \rm
Let $(A,\Delta)$ be a $^*$-algebraic quantum group. If $\omega=\varphi(\cdot\,a)$, then 
$$\widehat\psi(\omega^*\omega)=\varphi(a^*a).$$
\snl\bf Proof: \rm
By the formula in Lemma 2.1,we have
$$\widehat\psi(\omega^*\omega)=\omega^*(S^{-1}(a))=\omega(a^*)^-=\varphi(a^*a)^-=\varphi(a^*a).$$
\hfill $\square$
\einspr

This is indeed Plancherel's formula. Remember that in the case of Fourier series, it is written as
$$\sum_{n=-\infty}^{+\infty}|\widehat f(n)|^2=\frac{1}{2\pi}\int_0^{2\pi}|f(e^{i\theta})|^2\,d\theta.$$
As we mentoned already before (see Example 1.4), one needs to take proper normalizations of the dual Haar measures. In the case of Fourier analysis on $\Bbb R$, we have it under the form
$$\frac{1}{\sqrt{2\pi}}\int_{-\infty}^{+\infty}|\widehat f(p)|^2\,dp = 
                       \frac{1}{\sqrt{2\pi}}\int_{-\infty}^{+\infty}|f(t)|^2\,dt$$
where now 
$$\widehat f(p)=\frac{1}{\sqrt{2\pi}}\int_{-\infty}^{+\infty}f(t)e^{-itp}\,dt.$$
This case however does not fit in the purely algebraic theory. One needs to pass to the analytical theory of locally compact quantum groups (see further).

\nl
Among the $^*$-algebraic quantum groups, there is the distinquished class with positive integrals. It can be shown that for a positive left integral $\varphi$, the right integral $\varphi\circ S$ is again positive. The reason is that the scaling constant is forced to be $1$ (see e.g.\ [DC-VD]). This is not an obvious result. In [K-VD], it was already shown that a positive right integral exists if a positive left integral exists but it was not discovered yet that the scaling constant had to be $1$. However, from the positivity of $\varphi$, we easily get that also $\widehat\psi$ is positive. This is a consequence of the formula, just proven in Proposition 2.3 above. Also $\widehat\varphi$ is positive. In particular, the $^*$-algebraic quantum groups with positive integrals form a self-dual subfamily of the algebraic quantum groups. They are also locally compact quantum groups in the sense of Kustermans and Vaes, cf.\ [K-V1] and [K-V2].

\snl
Let me mention in passing that all of what I say in this talk can also be generalized to all locally compact quantum groups, but usually, these results are much more involved (and technically far more complicated). 

\nl
To finish this section, let us now consider the following well-known result in analysis: The {\it dual} of an abelian {\it compact group} is an abelian {\it discrete group}. What about the analogue in our setting?

\snl
We define an algebraic quantum group $(A,\Delta)$ of {\it compact type} if $A$ has an identity (i.e.\ when it is a Hopf algebra with integrals). This is motivated by the fact that the algebra $K(G)$ of continuous complex functions with compact support on a locally compact space $G$ has an identity if and only if $G$ is compact. On the other hand, an algebraic quantum group is called of {\it discrete type} if it has co-integrals. The use of this name can also be justified with reference to analysis, but this is less obvious. Just observe that the multiplier Hopf algebra of complex functions with finite support on a discrete group is indeed of discrete type.

\snl
Then we have the following result.

\inspr{2.4} Theorem \rm
If $A$ is of compact type, then $\widehat A$ is of discrete type. Similarly, if $A$ is of discrete type, then $\widehat A$ is of compact type. If $A$ is both of discrete and of compact type, $A$ has to be finite-dimensional.
\snl\bf Proof: \rm
If $1\in A$, then $\varphi\in \widehat A$ and by the left invariance of $\varphi$, we have
$$\omega\varphi=\omega(1)\varphi=\varepsilon(\omega)\varphi$$
so that $\varphi$ is a left co-integral in $\widehat A$.
\snl
If on the other hand, $h\in A$ and $ah=\varepsilon(a)h$ for all $a\in A$, then $\varphi(xh)=\varepsilon(x)\varphi(h)$ for all $x\in A$. Because $\varphi$ is faithful, we must have $\varphi(h)\neq 0$ as $h\neq 0$. Therefore, $\varepsilon\in \widehat A$ and so $\widehat A$ has a unit.
\snl
Moreover, in this case, we see that 
$$\align a\varphi(h)&=(\iota\otimes\varphi)((a\otimes 1)\Delta(h)) \\
                    &= (\iota\otimes\varphi)((1\otimes S^{-1}(a_{(2)}))\Delta(a_{(1)}h)) \\
                    &= (\iota\otimes\varphi) ((1\otimes S^{-1}(a))\Delta(h))
\endalign$$
for all $a\in A$ (we use $\iota$ for the identity map). So, if also $A$ has an identity, so that $\Delta(h)\in A\otimes A$, we see that any element $a$ of $A$ belongs to the left leg of $\Delta(h)$ and this is a finite-dimensional space. Hence $A$ is finite-dimensional.
\hfill $\square$
\einspr
We see that $A$ has to be finite-dimensional if it is both of compact and discrete type. This is the same as saying that, if we want both $A$ and $\widehat A$ to be Hopf algebras (with integrals), we must have that $A$ is finite-dimensional. This explains why it is important to work with multiplier Hopf algebras. 
\nl\nl

\bf 3. An example: The additive group $\Bbb Q_p$ of $p$-adic numbers. \rm
\nl
In this section, we will consider an example. It is well-known that the additive group $\Bbb Q_p$ of $p$-adic numbers (for any prime number $p$) is a locally compact group, it is non-compact and non-discrete but totally disconnected. Therefore, there is a natural algebraic quantum group associated with it. It is commutative and cocommutative, but it is neither of compact, nor of discrete type. For this reason, it can not be treated properly within the framework of Hopf algebras. Non-commutative and non-cocommutative examples can be constructed from the $p$-adic numbers using the basic construction methods. To keep things simple for this talk, I will only consider the basic case: the additive group $\Bbb Q_p$.
\nl
Recall that any $p$-adic number $x$ can be represented as a Laurent series
$$x=\sum_{j=-\infty}^{+\infty} x_j\,p^j$$
where $x_j\in \{0,1,2,\ldots,p-1\}$ and $x_j=0$ for $j$ small enough. The topology is determined by the metric $d$, defined by $d(x,y)=|x-y|$ where $|x|$ for any $x\in \Bbb Q_p$ is given by $p^{-m}$ with $m\in \Bbb Z$ satisfying $x_j=0$ if $j<m$ and $x_m\neq 0$. The $p$-adic numbers are defined as the completion of the rational numbers for this metric. Addition and multiplication on $\Bbb Q$ are continuous and therefore extend to $\Bbb Q_p$. Addition makes $\Bbb Q_p$ into an abelian group. As a topological space, it is non-discrete and non-compact but locally compact and totally disconnected. The $p$-adic integers $\Bbb Z_p$ are the numbers $x$ with $x_j=0$ when $j<0$. They form a compact open subgroup of $\Bbb Q_p$. In fact, the
subsets $p^m \Bbb Z_p$ with $m\in \Bbb Z$ form a basis for the topology around $0$ of compact open subgroups.  
\snl
For any totally disconnected locally compact group, we have the following result (see [L-VD2]).

\inspr{3.1} Proposition \rm Let $G$ be a totally disconnected locally compact group. Then the $^*$-algebra $A$ of continuous complex functions with compact support on $G$ is a $^*$-algebraic quantum group.
\einspr

The coproduct is the usual one: For $f\in A$, we have $\Delta(f)(r,s)=f(rs)$ whenever $r,s\in G$. Precisely because $G$ is totally disconnected, we have that functions like $\Delta(f)(1\otimes g)$ belong to the (algebraic) tensor product $A\otimes A$. The left integral is obtained by integration with respect to the left Haar measure on $G$. 
\snl
In the case $\Bbb Q_p$, the continuous complex functions with compact support are nothing else but functions $f$ with the property that $n,m$ in $\Bbb Z$ exist with $n<m$ and such that $f=0$ outside $p^n\Bbb Z_p$ and that $f$ is constant on $p^m\Bbb Z_p$.  
\snl
The Pontryagin dual of $\Bbb Q_p$ is again $\Bbb Q_p$. The bicharacter $\chi$, realizing this duality, is given by
$$\chi(x,y)=\exp(2\pi i xy),$$ 
where it is understood that 
$$\exp(2\pi i z)=\exp(2\pi i \sum_{j<0}z_j p^j)$$
if $z\in \Bbb Q_p$ is represented as $\sum_j z_j p^j$. If the Haar measure on $\Bbb Q_p$ is normalized so that the measure of the $p$-adic integers $\Bbb Z_p$ is $1$, then the dual Haar measure is the same.

\nl
We will use this example to illustrate the general result below (Proposition 3.3). First we need a definition. As before, we consider the case of a $^*$-algebraic quantum group $(A,\Delta)$ with {\it positive} integrals.

\inspr{3.2} Definition \rm
A non-zero element $h$ in $A$ is called a {\it group-like projection} if $h^2=h=h^*$ and 
$$\Delta(h)(1\otimes h)=h\otimes h.$$
\einspr

If $(A,\Delta)$ is the multiplier Hopf $^*$-algebra of complex functions with finite support on a (discrete) group $G$, then $h$ will be a group-like projection as in Definition 3.2 if and only if there is a finite subgroup $K$ of $G$ such that $h$ is the characteristic function of $K$ (i.e.\  the function with value $1$ on $K$ and $0$ everywhere else). In the case of $\Bbb Q_p$, we get examples of such group-like projections by considering the characteristic functions of the compact open subgroups.
\snl
Group-like projections in algebraic quantum groups are studied in [L-VD1]. They give rise to many objects (like the quantum analogue of the algebra of functions that are constant on left and right cosets associated with a subgroup). 
They also led to a theory of algebraic quantum hypergroups (see [De-VD1]). And it is expected that also a theory of  totally disconected locally compact quantum groups will result from this work. 
\snl
One property in this theory is related with the Fourier transform. Let us therefore consider it here.

\inspr{3.3} Proposition \rm
Let $h$ be a group-like projection in a $^*$-algebraic quantum group $(A,\Delta)$ with positive integrals. Then its Fourier transform $\widehat h$ is again a group-like projection in the dual $(\widehat A,\widehat \Delta)$ (provided we normalize the left integral $\varphi$ so that $\varphi(h)=1$).

\snl\bf Proof: \rm
For any $a\in A$ we have 
$$ \langle a,\widehat h\,\widehat h\rangle 
        =(\varphi\otimes\varphi)(\Delta(a)(h\otimes h))$$ 
and because also $\Delta(h)(h\otimes 1)=h\otimes h$, we get, using the left invariance of $\varphi$, that this is equal to
$$ (\varphi\otimes\varphi)(\Delta(ah)(h\otimes 1))=\varphi(h)\varphi(ah).$$
Because $\varphi(h)=\varphi(h^*h)$, we must have that $\varphi(h)\neq 0$. So, if we normalize $\varphi$ such that $\varphi(h)=1$, we obtain that $\widehat h^2=\widehat h$. We also have for any $a\in A$
$$\langle a, \widehat h^*\rangle=\langle S(a)^*,h\rangle^-=\varphi(S(a)^*h)^-=\varphi(hS(a)).$$
Because we have $S(h)=h$ and also $\varphi(S(ah)=\varphi(ah\delta)=\varphi(ah)$, we find that $\widehat h^*=\widehat h$.
Finally, for all $a,b\in A$ we get
$$\align \langle a\otimes b,\widehat \Delta(\widehat h)(1\otimes \widehat h)\rangle 
                      &=\langle (a\otimes 1)\Delta(b),\widehat h\otimes \widehat h \rangle \\
                      &=(\varphi\otimes\varphi)((a\otimes 1)\Delta(b)\Delta(h)(h\otimes 1)) \\
                      &=\varphi(ah)\varphi(bh)
\endalign$$
and this shows that $\widehat\Delta(\widehat h)(1\otimes \widehat h)=\widehat h\otimes \widehat h$.
\einspr
Let us illustrate this result in the case of the $p$-adic numbers. Consider for any $n\in \Bbb Z$ the compact open subgroup $K_n=p^n\Bbb Z_p$. Normalize the Haar measure so that the measure of the $p$-adic integers $\Bbb Z_p$ is $1$. Then, the measure of $K_n$ is $p^{-n}$. Denote by $h_n$ the characteristic function of $K_n$. It is a group-like projection in the associated algebraic quantum group because $K_n$ is a compact open subgroup of $\Bbb Q_p$. The Fourier transform $\widehat h_n$ of $h_n$ is given by
$$\align \widehat h_n(y) &=\int_{\Bbb Q_p} h_n(x) e^{-2\pi i xy}\, dx \\
           &= \int_{ p^n \Bbb Z_p} e^{-2\pi i xy}\, dx \\
           &= \frac{1}{p^n}\int_{\Bbb Z_p} e^{-2\pi i p^nxy}\, dx
\endalign$$
(where we are integrating over the Haar measure on $\Bbb Q_p$). 
We get $0$ except if $y\in p^{-n}\Bbb Z_p$ and then we get $p^{-n}$. Therefore $\widehat h_n=p^{-n}h_{-n}$.
\nl\nl

\bf 4. Conclusions and further research\rm
\nl
We have given the definition of the Fourier transform in the case of a multiplier Hopf algebra with integrals (in Section 1) and we have proven its basic and elementary properties (in Section 2). Multiplier Hopf $^*$-algebras with positive integrals are a special case of locally compact quantum groups. The Fourier transform and its properties can also be obtained for these locally compact quantum groups but the theory gets far more complicated and uses highly non-trivial operator algebra techniques. However, only here it gets its full strength including all aspects of abstract harmonic analysis.
\snl
On the other hand, it seems that the multiplier Hopf algebras with integrals, where everything remains purely algebraic, are sufficiently general, not only to function as a source of inspiration for the general theory, but also to provide already interesting special cases and examples.
\snl
We have illustrated this using the additive group $\Bbb Q_p$ of $p$-adic numbers. Using some of the basic constructions (like Majid's bicrossed product construction), it is possible to obtain more complicated examples. One could take e.g.\ a group like the $ax+b$-group, but build with the $p$-adic numbers in stead of the reals, and decompose it in order to get a matched pair of totally disconnected locally compact groups. Doing so, one can obtain interesting new and highly non-trivial examples of algebraic quantum groups. Some work along these lines has been done already, see e.g.\ [D-VD-W], but much more is possible.
\snl
In our last section, we have given the definition of a group-like projection in a $^*$-algebraic quantum group (with positive integrals). In [L-VD1], it is discovered that this relatively simple notion led to a variety of objects and results, involving the Fourier transform. Among other things, it was the starting point for the development of a theory of {\it algebraic quantum hypergroups}, see [De-VD1]. For such quantum hypergroups it turns out that a dual exists and that the duality obtained is generalizing in a nice way the duality we have already for the algebraic quantum groups. Again, more research is possible here (see e.g.\ [De-VD2]).
\snl
Another related topic is that of the study of all sorts of quantum subgroups of algebraic quantum groups and of the various associated objects (like the quantum version of the algebra of functions on cosets)
\snl
Finally, as we mentioned already, it is expected that all of this will eventually lead to a theory of totally disconnected locally compact quantum groups (examples of which have been discussed above). 

\nl\nl

\bf References \rm
\nl
\parindent 1.3 cm

\item{[DC-VD]} K.\ De Commer \& A.\ Van Daele: {\it Multiplier Hopf algebras embedded in C$^*$-algebraic quantum groups}. Preprint K.U.\ Leuven (2006), arxiv math.OA/0611872
\item{[De-VD1]} L.\ Delvaux \& A. Van Daele: {\it Algebraic quantum hypergroups}. Preprint University of Hasselt and K.U.\ Leuven (2006), arxiv math.RA/0606466
\item{[De-VD2]} L.\ Delvaux \& A. Van Daele: {\it A construction of algebraic quantum hypergroups by group-like idempotents and conditional expectations}. Preprint University of Hasselt and K.U.\ Leuven (2006).
\item{[D-VD-W]} L.\ Delvaux, A.\ Van Daele \& S.\ Wang: {\it Bicrossproduct of multiplier Hopf algebras}. Preprint University of Hasselt, K.U.\ Leuven \& Southeast University Nanjing (2007).
\item{[F]} G.B.\ Folland: {\it A course in abstract harmonic analysis}. Studies in advanced mathematics. CRC Press, Boca Raton, London (1995).
\item{[K-V1]} J.\ Kustermans \& S.\ Vaes: {\it A simple definition for locally compact quantum groups}. C.R.\ Acad.\ Sci.\ Paris S\'er I 328 (1999), 871--876.
\item{[K-V2]} J.\ Kustermans \& S.\ Vaes: {\it Locally compact quantum groups}. Ann.\ Sci.\ \'Ecole Norm.\ Sup. 33 (2000), 837--934.
\item{[K-VD]} J.\ Kustermans \& A.\ Van Daele: {\it C$^*$-algebraic quantum groups arising from algebraic quantum groups}. Int.\ J.\ Math.\ 8 (1997), 1067--1139.
\item{[L-VD1]} M.B.\ Landstad \& A.\ Van Daele: {\it Compact and discrete subgroups of algebraic quantum groups}. Preprint NTNU (Trondheim) and K.U.\ Leuven (2006), arxiv math.OA/0702525
\item{[L-VD2]} M.B.\ Landstad \& A.\ Van Daele: {\it Groups with compact open subgroups and multiplier Hopf $^*$-algebras}. Preprint NTNU (Trondheim) and K.U.\ Leuven (2006), arxiv math.OA/0701525
\item{[VD1]} A.\ Van Daele: {\it Multiplier Hopf algebras}. Trans.\ Am.\ Math.\ Soc.\ 342 (1994), 917--932.
\item{[VD2]} A.\ Van Daele:  {\it An algebraic framework for group duality}.  Adv.\ Math.\ 140 (1998), 323--366.
\item{[VD3]} A.\ Van Daele: {\it Locally compact quantum groups. A von Neumann algebra approach}. Preprint K.U.\ Leuven (2006), arxiv math.OA/0602212
\item{[VD-W]} A.\ Van Daele \& Shuanhong Wang: {\it The Larson-Sweedler theorem for multiplier Hopf algebras}. J.\ Alg.\ 296 (2006), 75--95.
\item{[VD-Z]} A.\ Van Daele \& Y.\ Zhang: {\it A survey on multiplier Hopf algebras} In 'Hopf algebras and Quantum Groups', eds. S.\ Caenepeel \& F.\ Van Oyestayen, Dekker, New York (1998), pp. 259--309.

\end

\end